\documentclass[12pt]{article}
\usepackage{amssymb}
\usepackage{graphicx}
\usepackage{subcaption}
\usepackage{amsfonts}
\usepackage{latexsym}
\usepackage{amsthm}
\usepackage{amsmath}

\usepackage{amssymb, amscd}

\usepackage{url}
\usepackage[T1]{fontenc}
\usepackage{color}
\usepackage{tabu}

\usepackage{bbm}

\usepackage{amsmath,amssymb,amscd,}
\usepackage[all]{xy}
\usepackage{color}

\numberwithin{equation}{section}

\newtheorem{conj}{Conjecture}[section]
\newtheorem{thm}[conj]{Theorem}
\newtheorem{rem}[conj]{Remark}

\newtheorem{defin}[conj]{Definition}
\newtheorem{prop}[conj]{Proposition}
\newtheorem{cor}[conj]{Corollary}
\newtheorem{lema}[conj]{Lemma}

\newtheorem{exam}[conj]{Example}

\begin{document}

\date{February 25, 2021}

\title{Optimal colored Tverberg theorems \\ for prime powers}


\author{{Du\v{s}ko Joji\'{c}} \\
 {\small Faculty of Science}\\[-2mm]
 {\small University of Banja Luka}
\and Gaiane Panina\\
{\small St. Petersburg State University} \\[-2mm]
{\small St. Petersburg Department of}\\[-2mm] {\small Steklov Mathematical Institute}
\and Rade T. \v{Z}ivaljevi\'{c}\\
{\small Mathematical Institute}\\[-2mm] {\small SASA,
  Belgrade}\\[-2mm]}

\maketitle
\begin{abstract}\noindent
The colored Tverberg theorem of Blagojevi\'{c}, Matschke, and Ziegler (Theorem \ref{A}) provides optimal bounds for the colored Tverberg problem, under the condition that the number of intersecting rainbow simplices $r=p$ is a prime number.

Our Theorem \ref{thm:main-p-power-mapping-red} extends this result to an optimal colored Tverberg theorem for multisets of colored points, which is valid for each prime power $r=p^k$, and includes   Theorem \ref{A} as a special case for $k=1$.
One of the principal new ideas is to replace the ambient simplex $\Delta^N$, used in the original Tverberg theorem, by an ``abridged simplex'' of smaller dimension, and to compensate for  this reduction by allowing  vertices to repeatedly appear a controlled number of times in different rainbow simplices.
Configuration spaces, used in the proof, are combinatorial pseudomanifolds which can be represented as  multiple chessboard complexes. Our main topological tool is the   Eilenberg-Krasnoselskii theory of degrees of equivariant maps for non-free actions.

A quite different generalization arises if we consider  colored classes  that are (approximately) two times  smaller than in the classical colored Tverberg theorem. Theorem \ref{ThmBalancedColored}, which unifies and extends some earlier results of this type, is based on the constraint method and uses the high connectivity of the configuration space.
\end{abstract}

\section{Introduction}\label{sec:intro}

The following result is known as {\em the topological Tverberg theorem},  \cite{BSS,  oz87, Vol96-1, M}.

\begin{thm}\label{thm:TT} {\rm (\cite{BSS,  oz87, Vol96-1})}
Let $r = p^k$ be a prime power, $d\geq 1$, and $N = (r-1)(d+1)$. Then
for every continuous map  $f : \Delta^N \rightarrow \mathbb{R}^d$,   defined on an  $N$-dimensional simplex, there exist disjoint faces $\Delta_1,\dots, \Delta_r$ of $\Delta^{N}$ such that
\[
 f(\Delta_1)\cap\dots\cap f(\Delta_r) \neq \emptyset \, .
 \]
\end{thm}

It is known that the condition on $r$ is  essential. Indeed, as demonstrated in \cite{bfz2},   if $r$ is not a prime power the topological Tverberg theorem fails if $d$ is sufficiently large.

\medskip
The following relative of Theorem \ref{thm:TT} is sometimes referred to as the {\em  Optimal colored Tverberg theorem} \cite{BMZ},  see also the review paper \cite{Z17} (and our Section \ref{sec:conf-spaces}) where it is classified as a {\em type C colored Tverberg theorem}.

\begin{thm} {\rm (\cite{BMZ})}\label{A}
Let $r\geq 2$ be a prime, $d\geq 1$, and $N:=(r-1)(d+1)$. Let $\Delta^N$ be an $N$-dimensional simplex with a partition (coloring) of its vertex set into $d+2$ parts,
\[
V  = [N+1] = C_0\sqcup \dots \sqcup C_d\sqcup C_{d+1}\, ,
 \]
 with $\vert C_i\vert = r-1$ for $i\leq d$ and $\vert C_{d+1}\vert =1$. Then for every continuous map $f : \Delta^N \rightarrow \mathbb{R}^d$, there are $r$ disjoint ``rainbow simplices''  $\Delta_1,\dots, \Delta_r$ of $\Delta^N$ satisfying
\[
f(\Delta_1) \cap \dots \cap f(\Delta_r) \neq\emptyset \,
\]
where by definition a face $\Delta$ of $\Delta^N$ is a rainbow simplex if and only if  $\vert \Delta\cap C_j\vert \leq 1$ for each $j=0,\dots, d+1$.
\end{thm}

Note that Theorem \ref{A} does not include Theorem \ref{thm:TT} as a special case. Indeed,  we need a stronger condition in the  colored Tverberg theorem,  where $r$ is a prime rather than a prime power. It remains an interesting question if this condition on $r$ can be relaxed. 

\medskip
In the following section we argue that a proper extension of Theorem \ref{A} may require multisets of colored points. This may not be an accident, as it appears to be dictated by the topology of a naturally  associated configuration space. For a motivating example, comparing the old and new results, see the end of the next section.

\medskip

\subsection*{Optimal colored Tverberg theorems  for multisets of points}
Our (first) main new result (Theorem \ref{thm:main-p-power-mapping-red}) is valid for each prime power $r=p^k$, and includes   Theorem \ref{A} as a special case for $k=1$.
One of the guiding ideas in the proof is to replace the simplex $\Delta^N$ (used in both Theorems \ref{thm:TT} and \ref{A}) by a simplex of smaller dimension, and to compensate this by allowing its vertices to appear a controlled number of times in different faces $\Delta_i$.

\begin{thm}\label{thm:main-p-power-mapping-red}
Let $r=p^k$ be a prime power,  $d\geq 1$, and $N:=k(p-1)(d+1)$.
Let $\Delta^N$ be an $N$-dimensional simplex whose vertices are  colored by $d+2$ colors,  meaning that there is a partition $V = C_0\sqcup C_1\sqcup\dots\sqcup C_{d}\sqcup C_{d+1}$ into $d+2$ monochromatic subsets.  We also assume that: \begin{enumerate}
                     \item[{\rm (1)}] Each of the sets $C_0,\dots, C_d$ has $(p-1)k$  vertices. The vertices in each $C_i$ are assigned multiplicities, as prescribed by the vector  $\mathbbm{L} = (1,p,\dots,p^{k-1})^{\times (p-1)}\in \mathbb{N}^{k(p-1)}$.
                     \item[{\rm (2)}] The (exceptional) color class $C_{d+1}$ contains a single vertex with multiplicity one.
                   \end{enumerate}

We claim that under these conditions for any continuous map $f:\Delta^N \rightarrow  \mathbb{R}^d$ there exist $r$ (not necessarily disjoint or even different) faces $\Delta_1, \dots, \Delta_r$ of $\Delta^N$ such that:

\begin{enumerate}
  \item[{\rm (A)}] $f(\Delta_1) \cap ... \cap f(\Delta_r)\neq \emptyset$.
  \item[{\rm (B)}] The number of occurrences of each vertex of $\Delta^N$ in all faces $\Delta_i$, does not exceed the prescribed multiplicity of that vertex.
  \item[{\rm (C)}]  All  faces $\Delta_i$  are multicolored or {\em rainbow simplices}, in the sense that their vertices have different colors,  $(\forall i)(\forall j) \, \vert Vert(\Delta_i)\cap C_j\vert\leq 1$.
  \end{enumerate}
\end{thm}

\medskip

  Theorem \ref{thm:main-p-power-mapping-red} has an  alternative formulation where the rainbow simplices $\Delta_i$ are faces of the original simplex $\Delta^{(r-1)(d+1)}$, rather than the faces of the abridged simplex $\Delta^{k(p-1)(d+1)}$  (used in Theorem \ref{thm:main-p-power-mapping-red}). While Theorem \ref{thm:main-p-power-mapping-red} is more intuitive and has a clear geometrical meaning, an advantage of Theorem \ref{thm:main-p-power-mapping} is that all simplices $\Delta_i$ are distinct and the result is closer in form to Theorems \ref{thm:TT} and \ref{A}.

  \begin{thm}\label{thm:main-p-power-mapping}
 Assuming that $r=p^k$ is a prime power and $r' = r-1$ let $$K = K_{r',r',\dots,r',1}\cong [r']\ast [r']\ast\dots\ast [r'] \ast [1] = [r']^{\ast(d+1)}\ast [1]$$ be a $(d+1)$-dimensional simplicial complex on a vertex set $V = V_0\sqcup V_1\sqcup\dots\sqcup V_{d+1}$,   divided into $(d+2)$ color classes where $\vert V_0\vert =\vert V_1\vert = \dots = \vert V_{d}\vert =r'$ and $\vert V_{d+1}\vert  =1$.
  Assume that $f : K\rightarrow \mathbb{R}^d$ is a continuous map which admits a factorization  $f = \widehat{f} \circ \alpha$ for some map $\widehat{f} : K_{k(p-1),k(p-1),\dots, k(p-1),1} \longrightarrow \mathbb{R}^d$ where
  \[
              \alpha : K_{r',r',\dots,r',1} \longrightarrow K_{k(p-1),k(p-1),\dots, k(p-1),1}
  \]
   is the simplicial map arising from a choice of a  $\mathbb{L}$-collapse map $\theta : [r'] \rightarrow [k(p-1)]$ (Definition \ref{def:collapse}). Then there exist $r$ pairwise vertex disjoint simplices ($r$ vertex-disjoint rainbow simplices) $\Delta_1,\dots, \Delta_r$ in $K$ such that
\begin{equation}\label{eqn:image-intersection}
          f(\Delta_1)\cap \dots \cap f(\Delta_r)\neq\emptyset \, .
\end{equation}
\end{thm}

The following extension of Theorem \ref{A} is the main result of \cite{BMZ}.

\begin{thm} {\rm (\cite{BMZ})}\label{AA}
Let $r\geq 2$ be a prime, $d\geq 1$, and $N:=(r-1)(d+1)$. Let $\Delta^N$ be an $N$-dimensional simplex with a partition  of its vertex set into $m+1$ parts,
\[
V  = [N+1] = C_0\sqcup \dots \sqcup C_{m}\, ,
 \]
with $\vert C_i\vert \leq r-1$ for $i=0,...,m$. Then for every continuous map $f : \Delta^N \rightarrow \mathbb{R}^d$, there are $r$ disjoint rainbow simplices  $\Delta_1,\dots, \Delta_r$ in $\Delta^N$ satisfying
\[
f(\Delta_1) \cap \dots \cap f(\Delta_r) \neq\emptyset \,
\]
\end{thm}

Theorem \ref{AA} is deduced from Theorem \ref{A} by a direct combinatorial argument.  By a similar reduction procedure we obtain a result (Theorem \ref{thm:cor-1} in Section \ref{Sec:Gen}) which extends Theorem \ref{thm:main-p-power-mapping-red},  and unifies Theorems \ref{A}, \ref{thm:main-p-power-mapping-red}, and \ref{AA}  in a single statement.

\medskip\noindent 
{\bf Example A:}
 The following example illustrates the optimality of Theorem \ref{thm:main-p-power-mapping-red} and illuminates its relationship with other statements of Tverberg type. Let $\mathcal{S} = \{A_i^{\alpha_i}\}_{i=0}^{d+1}$ be a {\em multiset} of points in $\mathbb{R}^d$ where  $\{A_i\}_{i=0}^d$ are the vertices of a non-degenerate simplex and $A_{d+1}$ its barycenter. Each of the points $A_i\, (i=0,\dots, d+1)$ is assigned a different color and a  multiplicity $\alpha_i \in \{1, r-1\}$, where $\alpha_i=1$ only for $i=d+1$.

If we allow a scattering of points, where each $\{A_i^{\alpha_i}\}$ is replaced by a collection $\{A_i^j\}_{j=1}^{\alpha_i}$ of (possibly distinct) points in $\mathbb{R}^d$ then:
\begin{enumerate}
  \item[{\rm (1)}] The classical (affine) Tverberg theorem guarantees the existence of a partition $\mathcal{S} = S_1\sqcup\dots \sqcup S_r$ such that  ${\rm Conv}(S_1)\cap \dots \cap {\rm Conv}(S_r)\neq \emptyset$;
  \item[{\rm (2)}] Theorem \ref{A} says that, under assumption that $r$ is a prime, such a partition exists with an additional property that each set $S_k$ has at most one point in each color;
  \item[{\rm (3)}] Theorem \ref{thm:main-p-power-mapping-red} claims that the conclusion of Theorem \ref{A} is true under assumption that $r$ is a prime power, if we allow only a partial scattering where in each color the points create clusters (partition) of the type $1^{p-1}p^{p-1}\ldots (p^{k-1})^{p-1}$. (The clusters themselves can be scattered in an arbitrary fashion.)
      \end{enumerate}

\subsection*{Colored Tverberg theorem for small color classes}

For comparison  we include a quite different  extension of the colored Tverberg theorem to prime powers.
We are able to avoid the factorization condition on the mapping $f$ (needed in Theorem \ref{thm:main-p-power-mapping}), at a price of increasing the number of color classes (and decreasing the number of vertices in each color class). 

\medskip
The following theorem can be interpreted as a prime power extension (relative) of the results from \cite[Section 9]{bfz1}.

\begin{thm}\label{ThmBalancedColored}
Let $r=p^\alpha$ be a prime power, $d\geq 1$, and $N:=(r-1)(d+2)$.
 Let $\Delta^N$ be an $N$-dimensional simplex with a partition (coloring) of its vertex set into $t$ color classes,
\[
V  = [N+1] = C_1\sqcup \dots \sqcup C_{t}\, ,
 \]
 with $\vert C_i\vert \leq q= \frac{r+1}{2}$ for all $i$. Then for every continuous map $f : \Delta^N \rightarrow \mathbb{R}^d$, there are $r$ disjoint ``rainbow simplices''  $\Delta_1,\dots, \Delta_r$ of $\Delta^N$ satisfying

 \begin{enumerate}
   \item[{\rm (1)}] \[
f(\Delta_1) \cap \dots \cap f(\Delta_r) \neq\emptyset \,
\]
where (as before)  a face $\Delta$ of $\Delta^N$ is a rainbow simplex if and only if  $\vert \Delta\cap C_j\vert \leq 1$ for each $j=1,\dots, t$.
   \item[{\rm (2)}] The dimension of each of the simplices $\Delta_i$ is at most $k$, where $k$ is defined from $rk+s=(r-1)d $,  $k > 0$, and $0 \leq s<r$.
   \item[{\rm (3)}] There are at most $s$ simplices whose dimensions are $k$.
 \end{enumerate}
 \end{thm}

\medskip
The organizations of the paper is the following.  The proof of Theorem \ref{thm:main-p-power-mapping-red} (and its equivalent form, Theorem \ref{thm:main-p-power-mapping}) is given in Section \ref{sec:proof-red}.
The role of chessboard complexes and their generalizations, as configuration spaces for theorems of Tverberg type,  is briefly  reviewed in \ref{sec:conf-spaces}.
In Section \ref{sec:conf-spaces} we formulate our main topological result of Borsuk-Ulam type (Theorem \ref{thm:main-p-power}), used in the proof of Theorem \ref{thm:main-p-power-mapping-red}, and its companion (Theorem \ref{thm:main-p-power-degree}),  about degrees of maps from multiple chessboard complexes. The proof of Theorems \ref{thm:main-p-power-degree} and \ref{thm:main-p-power} is postponed for Section \ref{sec:proof-degree}.

In Section \ref{sec:chess-multiple} we develop the theory of multiple chessboard complexes in the generality needed for applications in the proofs of Theorems \ref{thm:main-p-power-mapping-red}, \ref{thm:main-p-power-mapping} and \ref{thm:main-p-power-degree}. The focus is on multiple chessboard complexes which turn out to be pseudomanifolds (Sections \ref{sec:hierarchy} and \ref{sec:degree}).

In Section \ref{Sec:Gen} we formulate and prove extensions of Theorem \ref{thm:main-p-power-mapping-red} (Theorems \ref{thm:cor-1} and \ref{thm:main-p-power-generalization}). In Section \ref{sec:constraint} we show how the method of ``unavoidable complexes'' can be adapted for applications to problems about multisets of points.

In Section \ref{sec:comparison} we outline the proof of \cite[Theorem 2.1]{KB} (slightly extended to the case of pseudomanifolds), as  one of the central result illustrating the Eilenberg-Krasnoselskii comparison principle for degrees of equivariant maps, in the case of non-free group actions.

Finally, in Section \ref{sec:last} we give a proof of Theorem  \ref{ThmBalancedColored}.

\section{Chessboard complexes and equivariant maps}
\label{sec:conf-spaces}

The central role of {\em chessboard complexes}, as proper configuration spaces for {\em colored Tverberg problem and its relatives}, was recognized in \cite{ZV92} almost thirty years ago. To the present day these complexes remain, together with their generalizations (the multiple chessboard complexes)  in the focus of research in this area of geometric combinatorics.

\bigskip
Recall that the (standard) chessboard complex $\Delta_{p,q}$ is the complex of
all non-attacking placements of rooks in a $(p\times
q)$-chessboard (a placement is non-attacking if it is
not allowed to have more than one rook in the same row or in the
same column). More generally,  the multiple chessboard complex   $\Delta_{p,q}^{{\mathbbm{A}}, {\mathbbm{B}}}$   (see Section \ref{sec:chess-multiple}), where $\mathbbm{A}\in \mathbb{N}^q$ and $\mathbbm{B}\in \mathbb{N}^p$, arises if we allow more than one rook in each row (each column), where their precise number is determined by vectors $\mathbbm{A}$ and $\mathbbm{B}$.

\bigskip
Following \cite[Section 21.4]{Z17}, variants of the colored Tverberg problem are classified as {\em type A, type B, or type C}
depending on whether $m = d, m < d$ or $m > d$ (resp.), where $m + 1$ is the number
of colors and $d$ is the dimension of the ambient space.

\medskip
Central results in this area are the {\em topological type C colored Tverberg theorem} (Theorem 2.2 in \cite{BMZ}, reproduced here as Theorem \ref{A}) and the {\em topological type B colored Tverberg theorem} \cite{ZV92, VZ94}. Both of these results are obtained by applications of the {\em Configuration Space/Test Map scheme} involving chessboard complexes (see \cite{Z17}).

The associated {\em test maps} are respectively (\ref{eqn:nova})
 (in the type C case) and (\ref{eqn:nasa}) (for the type B result),
\begin{equation}\label{eqn:nova}
f : (\Delta_{r-1,r})^{\ast (d+1)}\ast [r]
\stackrel{\mathbb{Z}/r}{\longrightarrow} W_r^{\oplus (d+1)}
\end{equation}
\begin{equation}\label{eqn:nasa}
f : (\Delta_{2r-1,r})^{\ast
(k+1)}\stackrel{\mathbb{Z}/r}{\longrightarrow} W_r^{\oplus d}
\end{equation}
where $W_r\subset \mathbb{R}^r$  is the standard $(r-1)$-dimensional representation of the  group $S_r$.

Both theorems are consequences of the corresponding
Borsuk-Ulam-type statements claiming that in either case the
$S_r$-equivariant map $f$ must have a zero, where in (\ref{eqn:nova}) $r$ is a prime number, and  in  (\ref{eqn:nasa})) $r$ is a power of a prime.

\medskip
The following theorem  extends (\ref{eqn:nova}) and serves  as a basis for a new type C topological Tverberg theorem, which extends (in a natural way) the result of Blagojevi\' c,   Matschke and Ziegler to the prime power case.

\begin{thm}\label{thm:main-p-power}
 Let $G = (\mathbb{Z}_p)^k $ be a $p$-toral  group of order $r = p^k$. Let $\Delta_{k(p-1),p^k}^{\mathbbm{1};\mathbb{L}}$ be the multiple chessboard complex (based on a $k(p-1)\times p^k$ chessboard), where $ \mathbbm{1} = (1,\dots, 1)\in \mathbb{R}^{p^k}$ and $\mathbb{L} =  (1,p,\dots, p^{k-1})^{\times (p-1)}\in \mathbb{R}^{k(p-1)}$. Let $\partial\Delta_{[p^k]}\cong S^{p^k-2}$ be the boundary of a simplex with $p^k$ vertices.
 Then there does not exist a $G$-equivariant map
 \[
      f :   (\Delta_{k(p-1),p^k}^{\mathbbm{1};\mathbb{L}})^{\ast (d+1)}\ast [p^k] \longrightarrow (\partial\Delta_{[p^k]})^{\ast (d+1)}\cong (S^{p^k-2})^{\ast (d+1)} \cong S^{(p^k -1)(d+1)-1} \, .
 \]
\end{thm}

Theorem \ref{thm:main-p-power} is a consequence of the following theorem about degrees of equivariant maps.

\begin{thm}\label{thm:main-p-power-degree}
 Let $G = (\mathbb{Z}_p)^k $ be a $p$-toral  group of order $r = p^k$. Let $\Delta_{k(p-1),p^k}^{\mathbbm{1};\mathbb{L}}$ be the multiple chessboard complex (based on a $k(p-1)\times p^k$ chessboard), where $ \mathbbm{1} = (1,\dots, 1)\in \mathbb{R}^{p^k}$ and $\mathbb{L} =  (1,p,\dots, p^{k-1})^{\times (p-1)}\in \mathbb{R}^{k(p-1)}$. Let $\partial\Delta_{[p^k]}\cong S^{p^k-2}$ be the boundary of a simplex with $p^k$ vertices.
 Then  ${\rm deg}(f) \neq 0 \, (\mbox{{\rm mod} } p)$  for any $G$-equivariant map
 \[
      f :   (\Delta_{k(p-1),p^k}^{\mathbbm{1};\mathbb{L}})^{\ast (d+1)} \longrightarrow (\partial\Delta_{[p^k]})^{\ast (d+1)}\cong (S^{p^k-2})^{\ast (d+1)} \cong S^{(p^k -1)(d+1)-1} \, .
 \]
\end{thm}

\section{Chessboard pseudomanifolds}\label{sec:chess-multiple}

Following \cite{jvz, jvz2},  a multiple chessboard complex
$ \Delta_{m,n}^{\mathbb{K};\mathbb{L}}  = \Delta_{m,n}^{k_1,\dots,
k_n; l_1,\dots, l_m}$ is an abstract simplicial complex with vertices in  $[m]\times [n]$,   where the simplices have at most $k_i$ elements in the row $[m]\times \{i\}$
and at most $l_j$ elements in each column $\{j\}\times [n]$).

\medskip
We shall be mainly interested in complexes $ \Delta_{m,n}^{\mathbbm{1};\mathbb{L}}  = \Delta_{m,n}^{1,\dots, 1; l_1,\dots, l_m}$ where at most one rook is permitted in each of the rows of the chessboard $[m]\times [n]$.

\begin{prop}\label{prop:pseudomanifold}
 The multiple chessboard complex    $\Delta_{m,n}^{\mathbbm{1};\mathbb{L}}$ is a pseudomanifold if
\begin{equation}\label{eqn:uslov_za_mnogostrukost}
n = l_1+l_2+\dots+ l_m +1 \, .
\end{equation}
More precisely, the links of simplices of codimension $1$ and $2$ are spheres of dimensions $0$ and $1$, while in codimension $3$ may  appear both $2$-spheres and $2$-dimensional tori $T^2$.
\end{prop}

\medskip\noindent
{\bf Proof:} Let $S\in \Delta_{m,n}^{\mathbbm{1};\mathbb{L}}$ and let $s_i := \vert S\cap (\{i\}\times [n])\vert$.  The link ${\rm Link}(S)$ is clearly isomorphic to the multiple chessboard complex $\Delta_{m,n}^{\mathbbm{1};\mathbb{T}}$ where $\mathbb{T} = (t_1,\dots, t_m)$ and $t_i := l_i-s_i$. (Here we allow that $t_j = 0$ for some $j\in [m]$.) The proof is completed by an explicit description of all multiple chessboard complexes that arise as links of simplices in codimension  $\leq 3$.

\medskip
If ${\rm codim}(S) = 1$ then there exists $j_0$ such that $l_{j_0} = s_{j_0}+1$ and $l_j = s_j$ for each $j\neq j_0$. The condition (\ref{eqn:uslov_za_mnogostrukost}) guarantees  that $t_{j_0} =1$, which together with $t_j=0$ for $j\neq j_0$ implies ${\rm Link}(S)\cong  \Delta_{1,2}\cong S^0$.

 \medskip
If ${\rm codim}(S) = 2$ then there are two possibilities. Either (I) there exists $j_0$ such that $l_{j_0} = s_{j_0}+2$ and $l_j = s_j$ for each $j\neq j_0$, or (II) there exists $j_0\neq j_1$ such that both $l_{j_0} = s_{j_0}+1, l_{j_1} = s_{j_1}+1$ and  $l_j = s_j$ for each $j\neq j_0, j_1$. In the first case ${\rm Link}(S) \cong \partial \Delta_{[3]}\cong S^1$, while in the second ${\rm Link}(S) \cong \Delta_{2,3} \cong S^1$.

\medskip
If ${\rm codim}(S) = 3$ then the number of non-zero entries in the vector $\mathbb{T} = (t_1,\dots, t_m)$ is $1, 2$ or $3$. In the first case ${\rm Link}(S) \cong \partial \Delta_{[4]} \cong S^2$. In the second case ${\rm Link}(S) \cong \Delta_{2,4}^{\mathbbm{1};\mathbb{T}}$, where $\mathbb{T} = (2,1)$, hence $\Delta_{2,4}^{\mathbbm{1};\mathbb{T}}\cong S^2$ .

\medskip
Finally, in the third case ${\rm Link}(S) \cong \Delta_{3,4} \cong T^2$. \hfill $\square$

\subsection{Hierarchy of pseudomanifolds   $\Delta_{m,n}^{\mathbbm{1};\mathbb{L}}$}
\label{sec:hierarchy}

From here on we tacitly assume that the chessboard complex  $\Delta_{m,n}^{\mathbbm{1};\mathbb{L}}$ satisfies the condition $n = l_1+\ldots + l_m+1$. If $\mathbb{L} = \mathbbm{1}\in \mathbb{N}^m$ then $\Delta_{m,n}^{\mathbbm{1};\mathbb{L}} = \Delta_{n-1,n}$ is a standard chessboard complex \cite{BLVZ}, while in the case $m=1$ the complex $\Delta_{m,n}^{\mathbbm{1};\mathbb{L}} \cong \partial \Delta_{[n]}$
is the boundary sphere $\partial \Delta_{[n]} \cong S^{n-2}$ of the simplex $\Delta_{[n]}:= 2^{[n]}$.

\medskip
The pseudomanifolds   $\Delta_{m,n}^{\mathbbm{1};\mathbb{L}}$ form a poset category where the complexes  $\Delta_{n-1,n}$ and $\partial \Delta_{[n]}$ play the role of the initial and terminal object. The morphisms in this category are the {\em $\theta$-collapse maps} $\Omega_\theta$,  described in the following definition.

\begin{defin}\label{def:collapse}
Assuming  $m'\geq m$,  choose an epimorphism  $\theta :
[m']\rightarrow [m]$. Let $\widehat{\theta} : [m']\times
[n]\rightarrow [m]\times [n]$ be the associated  map of
chessboards where $\widehat{\theta}(i,j) = (\theta(i), j)$. We say that a sequence
$\mathbb{B} = (b_1,\ldots, b_{m})$ is obtained by a $\theta$-collapse from a sequence $\mathbb{A} = (a_1,\ldots, a_{m'})$ if $b_i = \sum_{\theta(j) =i}~a_j$. Define  $\Omega_\theta : \Delta_{m',n}^{\mathbbm{1};\mathbb{A}} \rightarrow
\Delta_{m,n}^{\mathbbm{1};\mathbb{B}}$ as the induced map of
multiple chessboard complexes where  $\Omega_\theta(S) :=
\widehat{\theta}(S)$, for each simplex $S\in  \Delta_{m,n}^{\mathbbm{1};\mathbb{A}})$. (Informally, the map $\Omega_\theta$ merges together some columns of  $\Delta_{m',n}^{\mathbbm{1};\mathbb{A}}$, as dictated by $\theta$.)
\end{defin}
The special cases $\Omega_\theta : \Delta_{n-1,n} \rightarrow
\Delta_{m,n}^{\mathbbm{1};\mathbb{L}}$ and $\Omega_\theta : \Delta_{m,n}^{\mathbbm{1};\mathbb{L}} \rightarrow \partial\Delta_{[n]}$ are of particular importance.   These maps are completely determined by the vector $\mathbb{L}$ and in this case both the corresponding $\Omega_\theta$ and the associated map $\theta : [n-1]\rightarrow [m]$ are referred to as $\mathbb{L}$-collapse maps.

\medskip

The group $S_n$, permuting the rows of the chessboard $[m]\times [n]$, acts on the multiple chessboard complex  $\Delta_{m,n}^{\mathbbm{1};\mathbbm{L}}$. The simplicial map  $\Omega_\theta : \Delta_{m',n}^{\mathbbm{1};\mathbb{L}'} \rightarrow \Delta_{m,n}^{\mathbbm{1};\mathbb{L}}$, associated to a collapse map $\theta : [m'] \rightarrow [m]$,  is clearly $S_n$-equivariant.

\subsection{Orientability of pseudomanifolds  $\Delta_{m,n}^{\mathbbm{1};\mathbb{L}}$ }

\begin{prop}\label{prop:orientation}
  The pseudomanifold   $\Delta_{m,n}^{\mathbbm{1};\mathbb{L}}$ is always orientable.  It has a fundamental class  $\tau \in H_d(\Delta_{m,n}^{\mathbbm{1};\mathbb{L}}; \mathbb{Z})\cong \mathbb{Z}$ where $d={\rm dim}(\Delta_{m,n}^{\mathbbm{1};\mathbb{L}})= n-2$. A permutation $g\in S_n$ reverses the orientation (changes the sign of $\tau$) if and only if  $g$ is odd.
\end{prop}

\medskip\noindent
{\bf Proof:}   Let  $\Omega_\theta :  \Delta_{m,n}^{\mathbbm{1};\mathbb{L}} \rightarrow \partial \Delta_{[n]}$
 be the collapse map associate to the constant map $\theta : [m] \rightarrow [1]$ (Definition \ref{def:collapse}). In other words $\Omega_\theta$ is the map induced by the projection $[m]\times [n] \rightarrow [1]\times [n]$ of chessboards, where a simplex $S\in \Delta_{m,n}^{\mathbbm{1};\mathbb{L}}$ is mapped to a simplex $S' \in \partial \Delta_{[n]}$  if and only if
\[
(\forall i\in [n])  \, [(\{i\}\times [m])\cap S \neq\emptyset \,  \Leftrightarrow \, i\in S']  \, .
\]
 Let $\widehat{S}$ be the  simplex $S \in \Delta_{m,n}^{\mathbbm{1};\mathbb{L}}$  oriented by listing its vertices in the increasing order of rows. Note that if $\Omega_\theta(S) = S'\in \partial \Delta_{[n]}$ then  $\Omega_\theta(\widehat{S})= \widehat{S'}$.

Choose an orientation $\mathcal{O}'$ on the sphere $\partial \Delta_{[n]}$ and use this orientation to define, via the collapse map $\Omega_\theta$, an orientation $\mathcal{O}$ on   $\Delta_{m,n}^{\mathbbm{1};\mathbb{L}}$. More explicitly, an ordered simplex $\widehat{S}$ is positively oriented with respect to $\mathcal{O}$ if and only if $\widehat{S'}$ is positively oriented with respect to the orientation $\mathcal{O}'$. It is not difficult to check that $\mathcal{O}$ is indeed and orientation on the pseudomanifold $\Delta_{m,n}^{\mathbbm{1};\mathbb{L}}$ which has all the properties  listed in Proposition \ref{prop:orientation}. \hfill $\square$

\begin{cor}
 As a consequence of Proposition \ref{prop:orientation} the $S_n$-pseudomanifolds $\Delta_{m,n}^{\mathbbm{1};\mathbb{L}}$ and $\Delta_{m',n}^{\mathbbm{1};\mathbb{L}'}$  are {\em concordant} in the sense that each $g\in S_n$ either changes the orientation of both of the complexes if none of them.
\end{cor}

\subsection{Degree of the collapse map $\Omega_\theta$}
\label{sec:degree}

In the following proposition we calculate the degree of the map
$\Omega_\theta$.

\begin{prop}\label{prop:degree} The degree of the map $\Omega_\theta : \Delta_{m',n}^{\mathbbm{1};\mathbb{A}} \rightarrow
\Delta_{m,n}^{\mathbbm{1};\mathbb{B}}$ is,
\begin{equation}\label{eqn:degree-1}
\mbox{\rm deg} (\Omega_\theta) =
\genfrac(){0pt}{0}{\mathbb{B}}{\mathbb{A}} =
\frac{b_1!\, b_2!\ldots b_{m}!}{a_1!\, a_2!\ldots a_{m'}!}.
\end{equation}
In the special case when $m=1$ we obtain that the degree of
the map $\Omega_\theta$ is the multinomial coefficient,
\begin{equation}\label{eqn:degree-2}
\mbox{\rm deg} (\Omega_\theta) = \frac{(a_1+ a_2+\ldots + a_{m'})!}{a_1!\,
a_2!\ldots a_{m'}!}
\end{equation}
and in the special case $a_1=a_2=\ldots=a_{m'}=1$  {\rm (\ref{eqn:degree-1})} reduces to the formula,
\begin{equation}\label{eqn:degree-3}
\mbox{\rm deg} (\Omega_\theta) =  b_1!\, b_2!\ldots b_{m}!.
\end{equation}
\end{prop}

\medskip\noindent
{\bf Proof:}
 Each simplicial map  $\Omega_\theta : \Delta_{m',n}^{\mathbbm{1};\mathbb{A}} \rightarrow \Delta_{m,n}^{\mathbbm{1};\mathbb{B}}$  is non-degenerate in the sense that it maps bijectively the top dimensional simplices of $\Delta_{m',n}^{\mathbbm{1};\mathbb{A}}$ to top dimensional simplices of  $\Delta_{m,n}^{\mathbbm{1};\mathbb{B}}$. Moreover, it is an orientation preserving map so in order to calculate the degree of $\Omega_\theta$ it is sufficient to calculate the cardinality of the preimage $\Omega_\theta^{-1}(c_0)$ of the barycenter $c_0$ of a chosen top dimensional simplex of $\Delta_{m,n}^{\mathbbm{1};\mathbb{B}}$.

\medskip
 Since the degree is multiplicative it is sufficient
to establish formula (\ref{eqn:degree-3}). A simple calculation shows that the cardinality of the set $\Omega_\theta^{-1}(c_0)$ is, in the case of a map  $\Omega_\theta : \Delta_{n-1,n}  \rightarrow \Delta_{m,n}^{\mathbbm{1};\mathbb{B}}$,  indeed given by the formula (\ref{eqn:degree-3}).
 \hfill $\square$

\section{Proof of Theorem \ref{thm:main-p-power-mapping-red}}
\label{sec:proof-red}

By convention $\Delta = \Delta_C$ is a simplex spanned by a set $C$, in particular $\Delta^N \cong \Delta_C$ where $C = Vert(\Delta^N) =  C_0\sqcup C_1\sqcup\dots\sqcup C_{d}\sqcup C_{d+1}$.

\medskip
Recall that a set $S\subset C$  (and the corresponding face $\Delta_S\subseteq \Delta_C$) is called a {\em rainbow set} (rainbow face) if $\vert S\cap C_i \vert \leq 1$ for all $ i=0,1,\dots, d+1$. It follows that the set of all rainbow simplices is a subcomplex of $\Delta_C$ which has a representation as a join of $0$-dimensional simplicial  complexes:
\begin{equation}\label{eqn:duga-1}
       \mathfrak{R} = \mathfrak{Rainbow} := C_0\ast C_1\ast \dots \ast C_d \ast C_{d+1} \subset \Delta_C \, .
\end{equation}
 By assumption $\vert C_i \vert = m:= k(p-1)$ for $i=0,1,\dots, d$ and $\vert C_{d+1}\vert = 1$, or more explicitly  $C_i = \{c^i_{\alpha, \beta}\} \,
  (0\leqslant \alpha \leqslant k-1;  1\leqslant \beta \leqslant  p-1)$ for all $0\leq i\leq d$, and $C_{d+1} = \{c_0\}$. Theorem \ref{thm:main-p-power-mapping-red} claims that for each continuous map $f : \Delta_C \rightarrow \mathbb{R}^d$ there exist rainbow faces $\Delta_1,\dots, \Delta_r \in \mathfrak{R}$ such that:
\begin{enumerate}
  \item[{\rm (1)}]   Vertex $c_0$ appears in at most one of the faces $\Delta_i$;
  \item[{\rm (2)}]   For all $i,\alpha,\beta$ the vertex $c^i_{\alpha,\beta}$ may appear in not more than $p^\alpha$ faces $\Delta_1, \dots, \Delta_r$;
  \item[{\rm (3)}]   $f(\Delta_1)\cap\dots\cap f(\Delta_r) \neq \emptyset$.
\end{enumerate}
An $r$-tuple $(\Delta_1, \dots, \Delta_r)$ of rainbow simplices is naturally associated to  the join $\Delta_1 \ast  \dots \ast  \Delta_r\in \mathfrak{R}^{\ast r}$.
Our immediate objective is to identify the subcomplex    $\mathfrak{R}^{\ast r}_\mathbb{L}\subset  \mathfrak{R}^{\ast r}$ which collects all $r$-tuples $(\Delta_1, \dots, \Delta_r)$ satisfying conditions (1) and (2).

\medskip
By assumption $\Delta_{i, \nu} : = \Delta_i\cap C_\nu$ is either empty or a singleton,  for each rainbow simplex $\Delta_i$. A moment's reflection reveals that the union $\cup \{\Delta_{i, \nu}\}_{i = 1}^{r}$ is a simplex in
$ \Delta_{k(p-1),p^k}^{\mathbbm{1};\mathbb{L}}$, for $0\leq\nu\leq d$ and a  simplex in     $[r] =  [p^k]$ if $\nu = d+1$. It immediately follows that
\[
\mathfrak{R}^{\ast r}_\mathbb{L}  \cong  (\Delta_{k(p-1),p^k}^{\mathbbm{1};\mathbb{L}})^{\ast (d+1)}\ast [r] \, .
\]

Let $\hat{f} : \mathfrak{R} \rightarrow \mathbb{R}^d$ be the restriction of the map $f : \Delta_C \rightarrow \mathbb{R}^d$.  The corresponding map defined on the $r$-tuples of rainbow simplices, satisfying conditions (1) and (2) is the map
\[
   \hat{F} : (\Delta_{k(p-1),p^k}^{\mathbbm{1};\mathbb{L}})^{\ast (d+1)}\ast [r] \longrightarrow (\mathbb{R}^d)^{\ast r} \, .
\]
By composing with the projection $ (\mathbb{R}^d)^{\ast r} \rightarrow  (\mathbb{R}^d)^{\ast r}/D$ (where $D\cong \mathbb{R}^d$ is the diagonal) and the embedding $(\mathbb{R}^d)^{\ast r}/D \hookrightarrow  (W_r)^{\oplus (d+1)}$, where $W_r\cong \mathbb{R}^r/\mathbb{R}$  is the standard $(r-1)$-dimensional representation of $S_r$, we obtain a map
\[
   \breve{F} : (\Delta_{k(p-1),p^k}^{\mathbbm{1};\mathbb{L}})^{\ast (d+1)}\ast [r] \longrightarrow (W_r)^{\oplus (d+1)}
\]
which has a zero in a simplex $(\Delta_1, \dots, \Delta_r)$ if and only if $f(\Delta_1)\cap \dots\cap f(\Delta_r)\neq \emptyset$. Since the sphere $S((W_r)^{\oplus (d+1)})\cong (S(W_r))^{\ast (d+1)}$ is equivariantly homeomorphic to $(\partial\Delta_{[r]})^{\ast (d+1)}$ a zero exists by Theorem \ref{thm:main-p-power}, which concludes the proof of Theorem \ref{thm:main-p-power-mapping-red}. \hfill $\square$

\section{Proof of Theorem \ref{thm:main-p-power-degree}}
\label{sec:proof-degree}

We are supposed to show that the degree ${\rm deg}(f)$ of each $G$-equivariant map
\begin{equation}\label{eqn:p-power-degree}
      f :   (\Delta_{k(p-1),p^k}^{\mathbbm{1};\mathbb{L}})^{\ast (d+1)} \longrightarrow (\partial\Delta_{[p^k]})^{\ast (d+1)}\cong (S^{p^k-2})^{\ast (d+1)} \cong S^{(p^k -1)(d+1)-1}
 \end{equation}
where $G = (\mathbb{Z}_p)^k$ is a $p$-toral group, is non-zero modulo $p$.
Following the {\em Comparison principle for equivariant maps} (Section \ref{sec:comparison})  we should:

\begin{enumerate}
  \item[{\rm (A)}]  Exhibit a particular map (\ref{eqn:p-power-degree}) such that ${\rm deg}(f) \neq 0$ modulo $p$;
  \item[{\rm (B)}]   Check if the conditions of Theorem \ref{thm:K-B-2.1-sing} are satisfied.
\end{enumerate}

The following proposition provides the needed example for the first part of the proof.

\begin{prop}\label{prop:map-example}
The $\theta$-collapse map
\begin{equation}\label{eqn:map-example}
\Omega_\theta : \Delta_{k(p-1),p^k}^{\mathbbm{1};\mathbb{L}} \longrightarrow \partial\Delta_{[p^k]} \cong \Delta_{1,r}^{1,r'},
\end{equation}
where $\theta: [k(p-1)] \rightarrow [1]$ is a constant map, has a non-zero degree modulo $p$.
\end{prop}

\medskip\noindent
{\bf Proof:} We calculate the degree of the map (\ref{eqn:map-example}) by applying the  formula (\ref{eqn:degree-2}). Recall that  $\mathbb{L} =  (1,p,\dots, p^{k-1})^{\times (p-1)}\in \mathbb{R}^{k(p-1)}$ so in this case
\begin{equation}\label{eqn:degree-4}
\mbox{\rm deg} (\Omega_\theta) = \frac{(p^k-1)!}{[(p^{k-1})!\, (p^{k-2})! \ldots p!\, 1!]^{p-1}} \, .
\end{equation}
The well-known formula for the highest power of $p$ dividing $m!$ is
\[
       {\rm ord}_p(m!)  = \left\lfloor \frac{m}{p}\right\rfloor + \left\lfloor \frac{m}{p^2}\right\rfloor + \dots  \, .
\]
By applying this formula we obtain
\[
     {\rm ord}_p((p^k-1)!) =  {\rm ord}_p((p^k)!) - k = p^{k-1} + p^{k-2} +\dots + 1 - k \]
and by applying the same formula to the denominator of (\ref{eqn:degree-4}) we obtain exactly the same quantity.   \hfill $\square$

\medskip
In light of Proposition \ref{prop:map-example} the map
\begin{equation}\label{eqn:collapse-(d+1)}
       (\Omega_\theta)^{\ast (d+1)}:   (\Delta_{k(p-1),p^k}^{\mathbbm{1};\mathbb{L}})^{\ast (d+1)} \longrightarrow (\partial\Delta_{[p^k]})^{\ast (d+1)}\cong (S^{p^k-2})^{\ast (d+1)} \cong S^{(p^k -1)(d+1)-1}
 \end{equation}
has a non-zero degree ${\rm deg}((\Omega_\theta)^{\ast (d+1)}) = ({\rm deg}(\Omega_\theta))^{ d+1}$ modulo $p$, which completes part (A) of the proof.

\medskip
For the part (B) of the proof of Theorem \ref{thm:main-p-power-degree} note (Section \ref{sec:comparison}) that it is sufficient to check the inequality (\ref{eqn:sufficient}). Let us begin with the observation that $\partial\Delta_{[r]}\, (r=p^k)$ is $S_r$-equivariantly homeomorphic to the unit sphere $S(W_r)$ in the standard $S_r$-representation  $W_r := \{x\in \mathbb{R}^r \mid \, x_1+\dots+ x_r = 0\}$.
As a consequence, for each subgroup $H\subseteq S_r$ the corresponding fixed point set $\partial\Delta_{[r]}^H \cong S(W_r)^H = S(W_r^H)$ is also a sphere.

\medskip
The action of $H$ decomposes $[r]$ into orbits $[r] = O_1\sqcup \dots \sqcup O_t$. From here easily follows a combinatorial description of the fixed point set $\partial\Delta_{[r]}^H$. A point $x\in \partial\Delta_{[r]}$, with barycentric coordinates $\{\lambda_i\}_{i=1}^r$, is fixed by $H$ if and only if the barycentric coordinates are constant in each of the orbits. Summarising, $\partial\Delta_{[r]}^H$ is precisely the boundary of the simplex  with vertices $\{o_i\}_{i=1}^t$, where $o_i$ is the barycenter of the face $\Delta_{O_i}\subset \Delta_{[r]}$.

  \medskip
 Let $\Delta_{m,r}^{\mathbbm{1};\mathbb{L}}$ be a multiple chessboard complex, where $\mathbbm{L} = (l_1,\dots, l_m)$ and $m = k(p-1)$.  It is not difficult to see that the barycenter $b_{i,j}$ of  (geometric realization of) the simplex $\{i\}\times O_j$ is in the fixed point set $(\Delta_{m,r}^{\mathbbm{1};\mathbb{L}})^H$ if and only if $\vert O_j\vert \leq l_i $.

 \medskip
  More generally,  a point $x$ is in $(\Delta_{m,r}^{\mathbbm{1};\mathbb{L}})^H$ if and only if it can be expressed as a convex combination
   \[
   x = \sum_{(i,j)\in S} \lambda_{i,j} b_{i,j}
   \]
  where $S$ is a subset of $[m]\times [t]$  satisfying
\begin{enumerate}
  \item[{\rm (1)}]   If $(i,j), (i',j)\in S$ then $i=i'$;
  \item[{\rm (2)}]  $(\forall i\in [m]) \, \sum \{\vert O_j\vert   \mid \, (i,j)\in S\}\leq l_i$.
\end{enumerate}
The $\theta$-collapse map $\Omega_\theta$, where $\theta: [m]\rightarrow [1]$ is the constant  map, maps $(\Delta_{m,r}^{\mathbbm{1};\mathbb{L}})^H$ to  $\partial\Delta_{[r]}^H$. Moreover $\Omega_\theta(b_{i,j}) = o_j$ and,  in light of (1) and (2), the simplex with vertices $\{b_{i,j}\}_{(i,j)\in S}$ is mapped bijectively to a face of  $\partial\Delta_{[r]}^H$. The following inequality is an  immediate consequence,
\[
     {\rm dim}((\Delta_{m,r}^{\mathbbm{1};\mathbb{L}})^H) \leq {\rm dim}(\partial\Delta_{[r]}^H) \, .
\]
From here and (\ref{eqn:collapse-(d+1)}) we obtain the inequality
  \[
        {\rm dim}  [(\Delta_{k(p-1),p^k}^{\mathbbm{1};\mathbb{L}})^{\ast (d+1)}]^H \leq {\rm dim}[(\partial\Delta_{[p^k]})^{\ast (d+1)}]^H
 \]
which finishes the proof of part (B) and concludes the proof of the theorem.
\hfill $\square$

\begin{rem}{\rm
  It follows from the part (B) of the proof of Theorem \ref{thm:main-p-power-degree} that $\Delta = (\Delta_{m,r}^{\mathbbm{1};\mathbb{L}})^H$ is also a ``chessboard complex''. Indeed,  $S\subseteq [m]\times [t]$ is a simplex in $\Delta$ if and only if $S$ has at most one rook in each row $[m]\times\{j\}$ and the total weight of the set  $(\{i\}\times  [t])\cap S$ is at most $l_i$, where the weight of each element $(i,j)$ is $\vert O_j\vert$.

  It follows that $\Delta$ can be classified as a complex of the type $\Delta_{m,t}^{\mathbbm{1}, \mathcal{L}}$ (cf. \cite[Definition 2.3]{jvz})), where $\mathcal{L}$ is family of threshold (simplicial) complexes.
  }
\end{rem}

\subsection{Proof of Theorem \ref{thm:main-p-power}}

If such a map $f$ exists, then there is a commutative diagram
  \begin{equation}\label{CD:3-to-2}
\begin{CD}
 (\Delta_{k(p-1),p^k}^{\mathbbm{1};\mathbb{L}})^{\ast (d+1)}\ast [p^k] @>{f}>> (\partial\Delta_{[p^k]})^{\ast (d+1)}\\
@A\pi AA @A\cong AA\\
 (\Delta_{k(p-1),p^k}^{\mathbbm{1};\mathbb{L}})^{\ast (d+1)} @>\widehat{f}>> (\partial\Delta_{[p^k]})^{\ast (d+1)}
\end{CD}
\end{equation}
where $\pi$ is the inclusion map and $\widehat{f} = f\circ\pi$. The map $\pi$ is homotopic to a constant map (since ${\rm Image}(\pi)$ is contained in ${\rm Star}(v)$ for any vertex $v\in [p^k]$). It follows  that ${\rm deg}(f) = 0$ which contradicts Theorem \ref{thm:main-p-power-degree}. \qed

\section{Generalizations by the method of constraints } \label{Sec:Gen}

In this section we prove Theorems \ref{thm:cor-1} and  \ref{thm:main-p-power-generalization} as extensions and relatives of Theorem \ref{thm:main-p-power-mapping-red}. Theorem \ref{thm:cor-1} unifies
the optimal colored Tverberg theorem (Theorem \ref{AA}), and the (primary) colored Tverberg theorem for multisets of points (Theorem \ref{thm:main-p-power-mapping-red}).

  \medskip
  All these results (in agreement with \cite{Z17}) can be  classified as  Type C colored Tverberg theorems (characterized by the condition that  the number of colors is at least $d+2$, where $d$ is the dimension of the ambient euclidean space).

  \medskip
 Theorems \ref{thm:cor-1} and  \ref{thm:main-p-power-generalization}  are deduced from Theorem \ref{thm:main-p-power-mapping-red} by a combinatorial reduction procedure closely related to the  method of ``constraining functions'' and ``unavoidable complexes'' \cite[Sections 3 and 4]{bfz1}, see also \cite{bestiary}, \cite{jmvz} and the review paper \cite{Z17} for more information. Here we show in Section \ref{sec:constraint} how  this method can be modified and extended to yield results about multisets of points.

\subsection{Unavoidable complexes }
\label{sec:constraint}

A multiset with vertices in $V$ is a pair $\mathbbm{V} = (V,m)$ where $m : V\rightarrow \mathbb{N}$  is a function assigning non-negative multiplicities to elements of $V$. If $V  = \{v_1,\dots, v_s\}$ and $m(v_i) = m_i$ we usually use the notation $\mathbbm{V} = \{v_1^{m_1},\dots, v_s^{m_s}\}$. We have introduced $\mathbbm{L}$ (in Theorem \ref{thm:main-p-power-mapping-red}) as the vector $(1,p,\dots,p^{k-1})^{\times (p-1)}\in \mathbb{N}^{k(p-1)}$. With a mild abuse of language we use the same notation for the corresponding multiset. More precisely a multiset is of the type $\mathbbm{L}$ if $(1,p,\dots,p^{k-1})^{\times (p-1)}$ is the associated vector of multiplicities.

\begin{defin}\label{def:unavoidable}{\rm (Unavoidable complexes)}
Let $\mathbbm{V} = (V,m)$ be a multiset and $r$ a positive integer. A simplicial complex $K\subseteq 2^V$ is $(r, \mathbbm{V})$-unavoidable  if for each $\mathbbm{V}$-proper collection  $\{\Delta_i\}_{i=1}^r$ of (not necessarily distinct) subsets of $V$,  at least one of the subsets $\Delta_i$ is in $K$. By definition a collection $\{\Delta_i\}_{i=1}^r$ is $\mathbbm{V}$-proper if for each $v\in V$  the cardinality of the set $\{i \mid v\in \Delta_i\}$ is at most $m(v)$.
  \end{defin}

\begin{defin}\label{def:T-complexes}{\rm (Tverberg complexes)}
Let $\mathbbm{V} = (V,m)$ be a multiset. Assume that $r$ and $d$ are positive integers. A simplicial complex $K\subseteq 2^V$ is a Tverberg complex of the type $(r,d,\mathbbm{V})$ if for each continuous map $f : K\rightarrow \mathbb{R}^d$
there exists a $\mathbbm{V}$-proper collection $\{\Delta_i\}_{i=1}^r$ of simplices in  $K$ such that
\begin{equation}\label{eqn:presek-postoji}
f(\Delta_1)\cap\dots\cap f(\Delta_r) \neq\emptyset \, .
\end{equation}
  \end{defin}

\begin{exam}\label{exam:T-type}{\rm
In the notation of Theorem \ref{thm:main-p-power-mapping-red} $$K_{d} := C_0\ast C_1\ast\dots\ast C_{d+1} \cong [k(p-1)]^{\ast (d+1)}\ast [1]$$
is a Tverberg complex of the type $(r, d, \mathbbm{V})$ where $r = p^k$ is a prime power and $\mathbbm{V}$ is a disjoint union of $d+1$ copies of the multiset  (of the type) $\mathbbm{L}$ and a singleton $[1]$.
}
\end{exam}

Unavoidable complexes, originally introduced in \cite[Definition 4.1]{bfz1} as ``Tverberg unavoidable subcomplexes'', play the fundamental role in the ``constraint method'' \cite[Sections 3 and 4]{bfz1}.
Here we show how the constraint method can be extended to the case of multisets.

\begin{prop}\label{prop:key}
  Let $\mathbbm{V} = (V,m)$ be a multiset with vertices in $V$.   Assume $K\subseteq 2^V$ is a Tverberg complex of the type $(r,d+1,\mathbbm{V})$, where  $r$ and $d$ are positive integers. Let $L$ be a $(r, \mathbbm{V})$-unavoidable complex. Then $K\cap L$ is a Tverberg complex of the type $(r,d,\mathbbm{V})$.
\end{prop}

\medskip\noindent
{\bf Proof:}  We are supposed to prove, following Definition \ref{def:T-complexes}, that  for each continuous map $f : K\cap L\rightarrow \mathbb{R}^d$
there exists a $\mathbbm{V}$-proper collection $\{\Delta_i\}_{i=1}^r$ of simplices in  $K\cap L$ which satisfies the condition (\ref{eqn:presek-postoji}).

The first step is to include $f : K\cap L \rightarrow \mathbb{R}^d$ into a commutative diagram (\ref{eqn:CD-1}) where $e$ and $i$ are the inclusion maps.
\begin{equation}\label{eqn:CD-1}
\begin{CD}
K @>F>> \mathbb{R}^{d+1}\\
@AeAA @AiAA\\
K\cap L @>f>> \mathbb{R}^d
\end{CD}
\end{equation}
Let $\bar{f}$ be an extension ($\bar{f}\circ e = f$) of
the map $f$ to $K$. Suppose that $\rho :
K\rightarrow \mathbb{R}$ is the function $\rho(x) := {\rm dist}(x,
K\cap L)$, measuring the distance of the point $x\in K$ from $K\cap L$.
Let $F = (\bar{f}, \rho)  : K\rightarrow \mathbb{R}^{d+1}$.

\medskip
By assumption  $K\subseteq 2^V$ is a Tverberg complex of the type $(r,d+1,\mathbbm{V})$ so there exists a $\mathbbm{V}$-proper family
$\{\Delta_1,\dots, \Delta_r\}$ of faces of $K$ such that
\begin{equation}\label{eqn:presek-postoji-2}
F(\Delta_1)\cap\dots\cap F(\Delta_r) \neq\emptyset \, .
\end{equation}
More explicitly, there exist $x_i\in\Delta_i$ such that $F(x_i)=F(x_j)$
for each $i,j = 1,\ldots, r$.

\medskip By assumption the complex $L$ is $(r, \mathbbm{V})$-unavoidable, hence $\Delta_i\in L$ for some $i\in [r]$. This implies $\rho(x_i)=0$ and
in turn $\rho(x_j)=0$ for each $j=1,\ldots, r$.

\medskip
If $\Delta_i'$ is the minimal face of $K$ containing $x_i$ then $\Delta_i'\in
L$ for each $i=1,\ldots, r$ and $f(\Delta_1')\cap\ldots\cap
f(\Delta_r')\neq\emptyset$.   \qed

\medskip
The proof of Proposition \ref{prop:key} is modeled on the proof of \cite[Lemma 4.3]{bfz1}. The following result is an immediate corollary,  see \cite[Theorem 4.4]{bfz1}.

\begin{cor}\label{cor:key}
  Let $\mathbbm{V} = (V,m)$ be a multiset with vertices in $V$.  Assume $K\subseteq 2^V$ is a Tverberg complex of the type $(r,d+c,\mathbbm{V})$, where  $r,d$ and $c$ are positive integers. Let    $L_1, \dots, L_c$  be a family of $(r, \mathbbm{V})$-unavoidable complexes. Then $K\cap L_1\cap \dots \cap L_c$ is a Tverberg complex of the type $(r,d,\mathbbm{V})$.
\end{cor}

\medskip\noindent
{\bf Proof:}  By induction, relying on Proposition \ref{prop:key}, we prove that $K\cap L_1\cap \dots \cap L_j$ is a Tverberg complex of the type $(r,d+c-j,\mathbbm{V})$ for each $j=1,\dots, c$.   \qed

\medskip
In the following proposition we exhibit a class of $(r, \mathbbm{V})$-unavoidable subcomplexes of $\Delta_V$, suitable for application of Corollary \ref{cor:key}.
(It should be compared to the example (i) in \cite[Lemma 4.2]{bfz1}.)

\begin{prop}\label{prop:unavoidable-examples-1}
  Let $\mathbbm{V} = (V,m)$ be a multiset with vertices in $V$.  Let $S\subset V$ be a subset such that
  \begin{equation}\label{eqn:m-weight}
        m(S):=\sum_{v\in S} m(v)   \leq r-1 \, .
  \end{equation}
  Then the complex $\Delta_{V\setminus S} = \{A\in 2^V \mid A\cap S =\emptyset\}$ is  $(r, \mathbbm{V})$-unavoidable.
  \end{prop}
In words, the face  $\Delta_{V\setminus S}$ of  $\Delta_{V}$ is a $(r, \mathbbm{V})$-unavoidable subcomplex of $\Delta_V$ if the total $m$-weight of $S$ does not exceed $r-1$.

\subsection{Extensions and relatives of Theorem \ref{thm:main-p-power-mapping-red}}

Our primary example of a Tverberg complex of the type $(r, d+c, \mathbbm{V})$ is the complex
\begin{equation}\label{eqn:initial}
K_{d+c} := C'_0\ast C'_1\ast\dots\ast C'_{d+c+1} \cong [k(p-1)]^{\ast (d+c+1)}\ast [1]
\end{equation}
described in Example \ref{exam:T-type} (with a slight change of notation). In this case the ambient simplex is $\Delta_V\cong \Delta^M$ where $M=k(p-1)(d+c+1)$,  $r = p^k$ is a prime power, and $\mathbbm{V}$ is a disjoint union of $d+c+1$ copies of the multiset $\mathbbm{L}$ with multiplicities $(1,  p, \dots, p^{k-1})^{p-1}$ and a singleton $[1]$,
\begin{equation}\label{eqn:main-multiset}
\mathbbm{V} = \mathbbm{L}^{\oplus (d+c+1)} \oplus [1] \, .
\end{equation}

\medskip
Choose pairwise disjoint subsets $S_1,\dots, S_c$ of $V$   such that the $m$-weight (\ref{eqn:m-weight}) of each $S_i$ does not exceed $r-1$. Let $S:= S_1\cup\dots\cup S_c$ and let $\Delta^N := \Delta_{V\setminus S}$.

\medskip
By Corollary \ref{cor:key} the complex
\[
     K_{d+c}\cap \Delta_{V\setminus S_1}\cap\dots\cap \Delta_{V\setminus S_c} = K_{d+c} \cap \Delta_{V\setminus S}
\]
is a Tverberg complex of the type $(r, d, \mathbbm{V})$. This result, formally translated in the language of Theorem \ref{thm:main-p-power-mapping-red}, may be reformulated as follows.

 \begin{thm}\label{thm:cor-1}
 Let $r=p^k$ be a prime power, $c\geq 1, d\geq 1$  and $N:=k(p-1)(d+1)$.
Let $\Delta^N$ be an $N$-dimensional simplex whose vertices are  colored by $d+c+2$ colors. More explicitly there is a partition $U = C_0\sqcup C_1\sqcup\dots\sqcup C_{d+c+1}$ of vertices of $\Delta^N$ into $d+c+2$ monochromatic sets (some  of the sets $C_i$ are allowed to be empty).   Let us assume that $\hat{m} : U\rightarrow \mathbb{N}$ is a function assigning positive multiplicities to the vertices of $\Delta^N$, which turns $U$  into a multiset $\mathbbm{U} = (U, \hat{m})$. Assume that
                     $\mathbbm{U} = (U, \hat{m})$ can be enlarged to the multiset $\mathbbm{V} = (V, m) = \mathbbm{L}^{\oplus (d+c+1)} \oplus [1] $  where $U\subseteq V$ and $m: V\rightarrow \mathbb{N}$ is an extension of $\hat{m}$. The enlargement is performed by adding of $c$ (pairwise disjoint) sets $S_1,\dots, S_c$    such that the $m$-weight (\ref{eqn:m-weight}) of each $S_i$ does not exceed $r-1$.

Then for any continuous map $f:\Delta^N \rightarrow  \mathbb{R}^d$ there exist $r$  faces $\Delta_1, \dots, \Delta_r$ of $\Delta^{N}$ such that:

\begin{enumerate}
  \item[{\rm (A)}] $f(\Delta_1) \cap ... \cap f(\Delta_r)\neq \emptyset$.
  \item[{\rm (B)}] The number of occurrences of each vertex of $\Delta^N$ in all faces $\Delta_i$, does not exceed the prescribed multiplicity of that vertex.
  \item[{\rm (C)}]  All  faces $\Delta_i$  are {\em rainbow simplices}, in the sense that   $(\forall i)(\forall j) \, \vert Vert(\Delta_i)\cap C_j\vert\leq 1$.\qed
  \end{enumerate}
\end{thm}
 Theorem \ref{thm:cor-1} reduces to Theorem \ref{AA} in the case $k=1$.

 \medskip

 If we assume in Theorem \ref{thm:cor-1} that $C_j = \emptyset$ for $j> d+1$ we obtain the following  extension of  Theorem \ref{thm:main-p-power-mapping-red}.  More directly it can be obtained from the case $c=1$ of Corollary \ref{cor:key}.
  Note that the number of colors and the total number of vertices (counted with multiplicities) is the same as in Theorem \ref{thm:main-p-power-mapping-red}. The main difference is that the color classes are treated equally (there are no ``exceptional'' colors).

\begin{thm}\label{thm:main-p-power-generalization}
Let $r=p^k$ be a prime power,  $d\geq 1$, and $N = k(p-1)(d+1)$.
Let $\Delta = \Delta^N$ be a simplex whose vertices are  colored by $d+2$ colors,  meaning that there is a partition $V = C_0\sqcup C_1\sqcup\dots\sqcup C_{d}\sqcup C_{d+1}$ into $d+2$ monochromatic subsets.  Assume that:
\begin{enumerate}
        \item[{\rm (1)}] The vertices of $C_i$  are assigned multiplicities from  the set $\{1,p,\dots, p^{k-1}\}$ so that each multiplicity is assigned to not more than $p-1$ elements of $C_i$;
  \item[{\rm (2)}] The total sum of all multiplicities over all the vertices of $\Delta$  is  $(r-1)(d+1)+1$.
\end{enumerate}

Under these conditions for any continuous map $f:\Delta \rightarrow  \mathbb{R}^d$ there exist $r$  faces $\Delta_1, \dots, \Delta_r$ of $\Delta$ such that:

\begin{enumerate}
  \item[{\rm (A)}] $f(\Delta_1) \cap ... \cap f(\Delta_r)\neq \emptyset$.
  \item[{\rm (B)}] The number of occurrences of each vertex of $\Delta^N$ in all faces $\Delta_i$, does not exceed the prescribed multiplicity of that vertex.
  \item[{\rm (C)}]  All  faces $\Delta_i$  are  rainbow simplices.
  \end{enumerate}
\end{thm}

\section{Comparison principle for equivariant maps }
\label{sec:comparison}

The following theorem is proved in \cite{KB} (Theorem 2.1 in Section 2). Note that the condition that the $H_i$-fixed point sets $S^{H_i}$ are {\em locally $k$-connected} for $k\leq {\rm dim}(M^{H_i})-1$ is automatically satisfied if $S$ is a representation sphere. So in this case it is sufficient to show that the sphere $S^{H_i}$ is (globally) $({\rm dim}(M^{H_i})-1)$-connected which is equivalent to the condition
\begin{equation}\label{eqn:sufficient}
{\rm dim}(M^{H_i}) \leq {\rm dim}(S^{H_i})\, (i=1,\dots, m) \, .
\end{equation}

\begin{thm}\label{thm:K-B-2.1}
Let $G$ be a finite group acting on a compact topological manifold $M = M^n$ and on a sphere $S \cong S^n$ of the same dimension.
Let $N\subset M$ be a closed invariant subset and let $(H_1), (H_2), \dots, (H_m)$ be the orbit types in $M\setminus N$. Assume that the set $S^{H_i}$ is both globally and locally $k$-connected for all $k=0, 1,\dots, {\rm dim}(M^{H_i})-1$, where $i = 1,\dots, m$.
Then for every pair of $G$-equivariant maps $\Phi, \Psi : M\longrightarrow S$, which are equivariantly homotopic on $N$, there is the following relation
\begin{equation}\label{eqn:fundam-congruence}
 {\rm deg}(\Psi) \equiv {\rm deg}(\Phi)  \quad ({\rm mod}\, GCD\{\vert G/H_1\vert, \dots, \vert G/H_k\vert\}) \, .
\end{equation}
\end{thm}

The proof of the following extension of Theorem \ref{thm:K-B-2.1} to manifolds with singularities doesn't require new ideas. By a {\em singular topological manifold} we mean a topological manifold with a codimension $2$ singular set. In particular Theorem \ref{thm:K-B-2.1-sing} applies to pseudomanifolds  $\Delta_{m,n}^{\mathbbm{1};\mathbb{L}}$, introduced in Section \ref{sec:chess-multiple}.

\begin{thm}\label{thm:K-B-2.1-sing}
Let $G$ be a finite group acting on a compact ``singular topological manifold'' $M = M^n$ and on a sphere $S \cong S^n$ of the same dimension.
Let $N\subset M$ be a closed invariant subset and let $(H_1), (H_2), \dots, (H_m)$ be the orbit types in $M\setminus N$. Assume that the set $S^{H_i}$ is both globally and locally $k$-connected for all $k=0, 1,\dots, {\rm dim}(M^{H_i})-1$, where $i = 1,\dots, m$.
Then for every pair of $G$-equivariant maps $\Phi, \Psi : M\longrightarrow S$, which are equivariantly homotopic on $N$, there is the following relation
\begin{equation}\label{eqn:fundam-congruence-sing}
 {\rm deg}(\Psi) \equiv {\rm deg}(\Phi)  \quad ({\rm mod}\, GCD\{\vert G/H_1\vert, \dots, \vert G/H_k\vert\}) \, .
\end{equation}
\end{thm}

\medskip\noindent
{\bf Proof:}  Following into footsteps of the proof of Theorem \ref{thm:K-B-2.1} (see \cite[Theorem 2.1]{KB})  we define a $G$-equivariant map
\begin{equation}
         f_0 : (M\times \{0,1\}) \cup (N\times [0,1]) \longrightarrow B\setminus \{O\}
\end{equation}
where  $B = {\rm Cone}(S)$ is a cone over the sphere $S$ (with the apex $O$), $\Psi$ and $\Phi$ are  restrictions of $f_0$ on $M\times \{0\}$ (respectively $M\times \{1\}$) and the restriction of $f_0$ on $N\times [0,1]$ is a homotopy between $\Psi \vert_N$ and $\Phi\vert_N$.

\medskip
If $f : M\times [0,1] \rightarrow B$ is a $G$-equivariant extension of $f_0$ then (\cite[Lemma 2.1]{KB})  ${\rm deg}(f) = \pm ({\rm deg}(\Psi) - {\rm deg}(\Phi))$
and the relation (\ref{eqn:fundam-congruence-sing}) will follow if
\begin{equation}\label{eqn:a-congruence}
  {\rm deg}(f) = \sum_{i=1}^{m}  a_i \cdot\vert G/H_i\vert
\end{equation}
for some integers $a_i\in \mathbb{Z}$.

\medskip
The proof of the following lemma (\cite[Lemma 2.2]{KB})  is quite general, in particular it holds for ``singular topological manifolds''.

\begin{lema}\label{lema:4-lema}
  There exists a $G$-equivariant extension $f : M\times [0,1] \rightarrow B$ of the map $f_0$ satisfying the following conditions:

  \begin{enumerate}
    \item[{\rm (}$\alpha${\rm )}] \,  $K = f^{-1}(O) = \bigcup_{j=1}^m T_j$ where $T_u\cap T_v = \emptyset$ for $u\neq v$;
    \item[{\rm (}$\beta${\rm )}] \, $T_j = G(K_j)$ for a compact set $K_j$;
    \item[{\rm (}$\gamma${\rm )}] \, $K_j = H_j(K_j)$ is $H_j$-invariant;
    \item[{\rm (}$\delta${\rm )}] \, $g(K_j)\cap h(K_j)=\emptyset$ if $gh^{-1}\notin H_j \, (j=1,\dots, m) $.
  \end{enumerate}
\end{lema}
The proof of Theorem \ref{thm:K-B-2.1-sing} is completed as in \cite[Section 2.1.3]{KB} by observing that ``singular topological manifolds'' also have absolute and relative fundamental classes.

\medskip
More explicitly, if $F_j$ is the restriction $f$ to a sufficiently small neighborhood of $K_j$ then
\[
    {\rm deg}(f)  = \sum_{j=1}^{m} {\rm deg}(F_j) \, .
\]
By the same argument as in \cite{KB} we deduce from Lemma \ref{lema:4-lema} that ${\rm deg}(F_j) = a_j\cdot \vert  G/H_j \vert$ for some $a_j\in \mathbb{Z}$, and the relation (\ref{eqn:fundam-congruence-sing}) is an immediate consequence.  \hfill $\square$

\section{Proof of Theorem \ref{ThmBalancedColored}}\label{sec:last}

The following theorem has already been formulated in \cite{JPVZ-Fomenko}, in a less general form. Note that even if $r$ is a prime, it does not immediately follow from Theorem \ref{A}.

\begin{thm} \label{C}
Let $r=p^\alpha$ be a prime power, $d\geq 1$, and $N:=(r-1)(d+1)$. Let $\Delta^N$ be an $N$-dimensional simplex with a partition (coloring) of its vertex set into $t$ color classes,
\[
V  = [N+1] = C_1\sqcup \dots \sqcup C_{t}\, ,
 \]
 where $\vert C_i\vert \leq q= \frac{r+1}{2}$ for all $i$. Then for every continuous map $f : \Delta^N \rightarrow \mathbb{R}^d$, there are $r$ disjoint \textit{rainbow simplices}   $\Delta_1,\dots, \Delta_r$ of $\Delta^N$ satisfying
\[
f(\Delta_1) \cap \dots \cap f(\Delta_r) \neq\emptyset \, .
\]
\end{thm}

\medskip\noindent
{\bf Proof:} If $m_i:= \vert C_i\vert $ then (by assumption) $\sum_{i=1}^{k}m_i=N+1$.
The configuration space of $r$-tuples of rainbow simplices is the join  $\mathcal{C}= \Delta_{r,m_1}\ast...\ast \Delta_{r,m_t}$. It can be visualized as the chessboard complex associated with $t$ ``small'' chessboards of height $r$ positioned side by side. The condition is that in each of these small chessboard the rooks placement is non-taking.

By assumption $r\geq 2m_i-1$ so by \cite[Proposition 1]{ZV92} $\Delta_{r,m_i}$ is $m_i-2$-connected. Since $\sum_{i=1}^{t}(m_i-2)+2(t-1)= (r-1)(d+1)-1$, it follows that  $\mathcal{C}$ is $(N-1)$-connected.

The reduction based on the standard   \textit{configuration space/test map scheme} \cite{M, Z17} shows that if $f$ violates the statement of the theorem, then there exists a $(\mathbb{Z}_p)^\alpha$-equivariant map $$\mathcal{C} \longrightarrow S^{(d+1)(r-1)-1}\, .$$
This contradicts Volovikov's theorem \cite{Vol96-1}. \qed

\bigskip\noindent
{\bf Proof of Theorem \ref{ThmBalancedColored}}.  As in the proof of Theorem \ref{C} let $\mathcal{C}= \Delta_{r,m_1}\ast...\ast \Delta_{r,m_t}$ be the configuration space of all $r$-tuples $(\Delta_1,\dots, \Delta_r)$ of (pairwise vertex disjoint) rainbow simplices, where $m_i:= \vert C_i\vert $ for each $i=1,\dots, t$.

As before $\sum_{i=1}^{k}m_i=N+1$ and the only difference with the setting of Theorem \ref{C} is that the dimension of the simplex $\Delta^N$ is now $N=(r-1)(d+2)$ (rather than $N=(r-1)(d+1)$).

\medskip
Let $\Delta_N^{(\nu)} = (\Delta^N)^{(\nu)}$ be the $\nu$-dimensional skeleton of the simplex $\Delta^N$ and let
\[
  \Sigma = SymmDelJoin(\Delta_N^{(k)},...,\Delta_N^{(k)}, \Delta_N^{(k-1)},...,\Delta_N^{(k-1)})
\]
be the \emph{symmetrized deleted join}, introduced in \cite{jvz2} (Theorem 3.3) as the configuration space of all  (pairwise vertex disjoint) $r$-tuples satisfying the conditions (2) and (3) in Theorem \ref{ThmBalancedColored}.

\medskip Our goal is to show that for each continuous map $f : \Delta^N \rightarrow \mathbb{R}^d$ there exists an $r$-tuple $\Delta\in \mathcal{C}\cap\Sigma$ such that $f(\Delta_1) \cap \dots \cap f(\Delta_r) \neq\emptyset$.
Suppose (for contradiction) that $f$ is a counterexample, violating  the statement of the theorem. By the standard reduction  $f$ induces an equivariant map

$$F: \mathcal{C} \rightarrow W_r^{\oplus (d+1)}$$ which does not have a zero in $\Sigma$. By \cite[Lemma 2.10]{F2} there exists
an equivariant  map $$\Phi: (\Delta^N)^{*r}_\Delta \rightarrow W_r$$ which has a remarkable property
$ \Phi^{-1}(0) = \Sigma$. It immediately follows that the map
$$(F,\Phi): \mathcal{C} \rightarrow W_r^{\oplus (d+2)}$$
has no zeros, however this is in contradiction with Theorem \ref{C}. \qed

\bigskip\noindent
{\bf Acknowledgements:}
R. \v Zivaljevi\' c was supported by the Serbian Ministry of Education, Science and Technological Development through Mathematical Institute of the Serbian Academy of Sciences and Arts. Section 8  is supported by the Russian Science Foundation under grant  21-11-00040.

\end{document}